\documentclass[reqno,centertags, 12pt]{amsart}

\usepackage{amssymb,amsmath,amsfonts,amssymb} 
-\marginparwidth \oddsidemargin -4mm \evensidemargin
\oddsidemargin
\usepackage{latexsym}
\advance\hoffset by 5mm

\def\@abssec#1{\vspace{.05in}\footnotesize \parindent .2in
{\bf #1. }\ignorespaces}
\newtheorem{theorem}{Theorem}[section]

\newtheorem{lemma}[theorem]{Lemma}

\def \Rm {\mathbb R}
\def \Tm {\mathbb T}

\def \Zm {\mathbb Z}

\def\cP{\mathcal P}
\def\cC{\mathcal C}

\allowdisplaybreaks \numberwithin{equation}{section}

\title[Global Regularity for
Dispersive SQG]{Global Regularity for the Critical Dispersive
Dissipative Surface Quasi-Geostrophic Equation}
\author{Alexander Kiselev and Fedor Nazarov}
\thanks{Department of
Mathematics, University of Wisconsin, Madison, WI 53706, USA;
email: kiselev@math.wisc.edu, nazarov@math.wisc.edu}

\begin{document}

\begin{abstract}
We consider surface quasi-geostrophic equation with dispersive
forcing and critical dissipation. We prove global existence of
smooth solutions given sufficiently smooth initial data. This is
done using a maximum principle for the solutions involving
conservation of a certain family of moduli of continuity.
\end{abstract}

\maketitle

\section{Introduction}\label{intro}

In this paper, we study the following dispersive dissipative
surface quasi-geostrophic (SQG) equation:
\begin{equation}\label{sqgd1}
\theta_t = u \cdot \nabla \theta - (-\Delta)^{1/2}\theta + Au_2,
\,\,\,\theta(x,0)=\theta_0(x).
\end{equation}
Here $\theta$ is a scalar real-valued function, $A$ is an
amplitude parameter, and the velocity $u$ is given by
$u=(-R_2\theta,R_1\theta)$ with $R_1,R_2$ the usual Riesz
transforms. We will consider \eqref{sqgd1} on a torus $\Tm^2$ (or,
equivalently, on $\Rm^2$ with periodic initial data). The equation
\eqref{sqgd1} arises in geophysical fluid dynamics and can be
interpreted as the evolution of buoyancy on a 2D surface in the
presence of an ambient buoyancy gradient. Physically, the presence
of the background gradient gives rise to dispersive waves and
hence, the system supports both wave-like and turbulent motions
(see Held et al. \cite{HPGS} for additional physical insight into
\eqref{sqgd1} and Sukhatme and Smith \cite{LSJS} for its
interpretation as part of a broader family).

In the recent years, the SQG equation has been focus of intense
mathematical research, initiated by Constantin, Majda and Tabak
\cite{CMT}. The equation is physically motivated, and it is perhaps
the simplest equation of fluid dynamics for which the question of
global existence of smooth solutions is still poorly understood.
Global regularity for the SQG equation without dispersion is known
in the subcritical regime, when the dissipative term is
$(-\Delta)^\alpha,$ $\alpha \geq 1/2$. The subcritical case
$\alpha>1/2$ goes back to Resnick \cite{Res} while the critical case
$\alpha = 1/2$ was recently settled in \cite{KNV} and \cite{CV}. The
supercritical case $\alpha < 1/2$ in general remains open.

Mathematically, the key ingredient in regularity proofs in the
subcritical case is the maximum principle for
$\|\theta(x,t)\|_{L^\infty}$ (see e.g. \cite{Res}). In the critical
case, the crucial improvement comes from the stronger nonlocal
maximum principle for a certain modulus of continuity (\cite{KNV})
or DiGiorgi-type iterative estimates establishing H\"older
regularity of $\theta$ (\cite{CV}). In this work, we extend the
nonlocal maximum principle technique to the case where dispersion is
also present. The main result is
\begin{theorem}\label{main}
The dispersive critical surface quasi-geostrophic equation
\eqref{sqgd1} with smooth periodic initial data has a unique
global smooth solution.
\end{theorem}
\noindent \it Remarks. \rm 1. We note that in the case of stronger
dissipation $\alpha > 1/2$ the result also holds true and can be
proven in a standard way (once the control of
$\|\theta(x,t)\|_{L^\infty}$ is established - which is a part of
our proof and can be extended to the subcritical case in a straightforward manner).  \\
2. The smoothness assumption on the initial data can be relaxed to
$\theta_0 \in H^1.$ Indeed, local existence of the solution smooth
for $t>0$ starting from such initial data can be proven by
standard methods. The linear dispersive part does not present any
difficulty in this respect (see e.g. \cite{D} for the SQG case or
\cite{KNS} for an argument in the case of Burgers equation, which
can be easily adapted to our situation). Once $t>0,$ one can apply
the Theorem above to get global smooth solution. \\
3. The proof of the uniqueness in the setting of Theorem~\ref{main}
is also standard. \\
4. The key step in the proof is, like in \cite{KNV}, the
derivation of a uniform estimate on $\|\nabla \theta\|_{L^\infty}$
by using a family of moduli of continuity preserved by the
evolution. Once one has this estimate, the proof of global
existence of regular solution is achieved by well-known approach
of using local existence theorem and  differential inequalities
for the Sobolev norms of the solution. Thus, in what follows we
will focus on the essential issue of gaining control of $\|\nabla
\theta\|_{L^\infty}.$

\section{The Proof}\label{proofs}

Our first observation is that the $L^2$
norm of the solution over a single period cell is non-increasing.

\begin{lemma}\label{L2}
The $L^2$ norm of a smooth solution of \eqref{sqgd1} is
non-increasing.
\end{lemma}
\begin{proof}
Multiplying the equation by $\theta(x,t)$ and integrating we
obtain
\[ \frac12 \partial_t \|\theta\|_{L^2} = -\int_{\Tm^2}
\theta(-\Delta)^{1/2}\theta\,dx + A\int_{\Tm^2} \theta u_2\,dx. \]
The first term on the right hand side is negative, while the
second is, up to a constant factor, equal to
\[ \sum_{k \in \Zm^2, k \ne 0}
\frac{k_1}{|k|}|\hat{\theta}(k_1,k_2)|^2. \] The latter expression
is zero since $\theta$ is real-valued and so $\hat{\theta}$ is
even.
\end{proof}

Our next step is gaining control of the $L^\infty$ norm of the
solution of \eqref{sqgd1}. One can no longer claim it is
non-increasing as in the non-dispersive case (it isn't), but it
remains uniformly bounded.

\begin{lemma}\label{Linf}
There exists a constant $D = D(A, \theta_0)$ such that the
$L^\infty$ norm of a smooth solution $\theta(x,t)$ of
\eqref{sqgd1} satisfies
\[ \|\theta(x,t)\|_{L^\infty} \leq D \]
for all times while the solution remains smooth.
\end{lemma}
\noindent \it Remark. \rm Observe that in contrast to the
non-dispersive case, there is no $L^\infty$ norm maximal
principle: the $L^\infty$ norm can grow, and numerical
computations suggest it often does \cite{LSJS}. Instead, we just
have an upper bound on the $L^\infty$ norm.
\begin{proof}
Consider a point $x$ where $\theta(x,t)$ reaches its maximum, $M$
(the case of a minimum is similar). At the point of maximum, we
have
\[ \partial_t \theta(x,t) = -(-\Delta)^{1/2} \theta(x,t) +A R_1
\theta(x,t). \]
We will use the representation
\begin{equation}\label{dissip} -(-\Delta)^{1/2} \theta(x,t) =
\bigl.\frac{d}{dh}\cP_h * \theta\bigr|_{h=0}= \lim_{h \rightarrow
0} \frac{1}{h} \int_{\Rm^2}
\frac{h}{(|y|^2+h^2)^{3/2}}(\theta(x-y)-M)\,dy, \end{equation}
where $\cP_h$ is the usual Poisson kernel in $\Rm^2.$ Observe that
we can pass to the limit in \eqref{dissip} obtaining the kernel
$|y|^{-3}.$ On the other hand,
\begin{equation}\label{R1}
A R_1 \theta(x,t) = A \int_{\Rm^2}
\frac{f(\hat{y})}{|y|^2}\theta(x-y)\,dy,
\end{equation}
where $f$ is a smooth mean zero function on the unit circle and
$\hat{y}=y/|y|.$ The integral converges in the principal value
sense. Because of the mean zero property of $f$ we can replace
$\theta(x-y)$ in \eqref{R1} with $\theta(x-y) -M.$

Consider a ball of radius $\rho$ centered at zero $B_\rho,$ and the
portion of the integrals in \eqref{dissip}, \eqref{R1} corresponding
to that ball:
\begin{equation}\label{ie1}
\int_{B_\rho} \left( \frac{1}{|y|^3} + \frac{Af(\hat{y})}{|y|^2}
\right) (\theta(x-y)-M)\,dy.
\end{equation}
We can choose $\rho=\rho(A)$ sufficiently small independently of
$M$ so that \eqref{ie1} does not exceed
\begin{equation}\label{ie2} \frac12 \int_{B_\rho}
\frac{1}{|y|^3} (\theta(x-y)-M)\,dy.
\end{equation}
Let us denote by $m$ the Lebesgue measure on $\Tm^2.$  Since by
Lemma~\ref{L2}, $\|\theta(x,t)\|_{L^2(\Tm^2)} \leq
\|\theta_0\|_{L^2(\Tm^2)},$ we have that
\begin{equation}\label{meascon}
m\left(y \in \Tm^2 \left| |\theta(y,t)| \geq
\frac{M}{2}\right)\right.\leq
\frac{4\|\theta_0\|_{L^2(\Tm^2)}^2}{M^2}.
\end{equation}
Assume that $\rho$ is sufficiently small so that $B_\rho$ fits
into a single period cell. Since
 $|y|^{-3}$ is monotone decreasing, the expression in \eqref{ie2}
is maximal if points where $\theta(y,t)$ is large are concentrated
near $y=0.$ In particular, assuming that $M$ is sufficiently
large, we see from \eqref{meascon} that the expression in
\eqref{ie2} is less than or equal to
\[ -\frac12 \int_{B_\rho \setminus B_r} \frac{M}{2
|y|^3}\,dy, \] where $r =
2\pi^{-1/2}\|\theta_0\|^2_{L^2(\Tm^2)}M^{-1}.$ Therefore, the
expression in \eqref{ie2} does not exceed
\begin{equation}\label{ie3} -
\frac{\pi}{2} \int_r^\rho \frac{M}{|y|^2}d|y| \leq
-\frac{\pi^{3/2}}{4\|\theta_0\|_{L^2(\Tm^2)}}M^2+C(A)M.
\end{equation}
The integral over the complement of $B_\rho$ in \eqref{dissip} is
negative, so it remains to control \begin{equation}\label{litest15}
A\int_{\Rm^2 \setminus
B_\rho}\frac{f(\hat{y})}{|y|^2}(\theta(x-y)-M)\,dy. \end{equation}
Note that due to the mean zero property of $f,$ we can replace $M$
in \eqref{litest15} with $\overline{\theta},$ the mean value of
$\theta$ over a period cell. Then for any period cell $\cC$ lying
entirely in $\Rm^2 \setminus B_\rho$ with center at distance $L$
from the origin, we have
\[ \left| \int_{\cC} \frac{f(\hat{y})}{|y|^2}(\theta(x-y)-\overline{\theta})\,dy
\right| \leq CM {\rm max}_{\cC} \left|\nabla
\left(\frac{f(\hat{y})}{|y|^2}\right)\right| \leq CML^{-3}. \]
Adding up contributions of the different cells, we get the total
bound $CAM$ for \eqref{litest15}.
 Therefore, we have
\[ \partial_t \theta(x,t) \leq -C(\theta_0)M^2 +C(A)M, \]
which is negative provided that $M$ is large enough; define
$\tilde{D}(A,\theta_0)$ so that this is true if $M \geq
\tilde{D}(A,\theta_0).$ But then it is clear that $\theta(x,t)$
can never reach such value of $M$ unless $\|\theta_0\|_{L^\infty}$
was already larger - but in this case, $\|\theta\|_{L^\infty}$
will decay until reaching at least $\tilde{D}(A,\theta_0).$
Setting $D(A,\theta_0) = {\rm
max}\{\tilde{D}(A,\theta_0),\|\theta_0\|_{L^\infty}\},$ we obtain
the result of the Lemma.
\end{proof}

Now we introduce a family of moduli of continuity. This is the
same family that was considered in \cite{KNV} in the case of the
critical SQG. Namely, let $\omega(\xi)$ be continuous and defined
by
\begin{eqnarray}\label{modcon1}
\omega(\xi) = \xi - \xi^{3/2}, \,\,\,\, & 0 \leq \xi \leq \delta; \\
\nonumber \omega'(\xi) = \frac{\gamma}{\xi(4+\log(\xi/\delta))},
\,\,\,\, & \xi \geq \delta,
\end{eqnarray}
and set $\omega_B(\xi) = \omega(B\xi).$ Here $0<\gamma<\delta$ are
certain constants defined in \cite{KNV}; the modulus of continuity
$\omega$ is increasing, concave and differentiable at every point
except $\xi =\delta.$

We will need the following lemma from \cite{KNV}:
\begin{lemma}\label{somod} If the function $\theta$ has
modulus of continuity $\omega_B$, then $u=(-R_2\theta, R_1\theta)$
has modulus of continuity
\[ \Omega_B(\xi)=C\left(\int_0^\xi \frac{\omega_B(\eta)}{\eta}\,d\eta+
\xi\int_\xi^\infty \frac{\omega_B(\eta)}{\eta^2}\,d\eta\right) \]
with some universal constant $C>0$.
\end{lemma}
Observe that by a simple change of coordinates and definition of
$\omega_B,$ $\Omega_B(\xi)=\Omega(B\xi).$

Our next lemma can be proven exactly as in \cite{KNV}, using that
$\omega''(0)=-\infty:$
\begin{lemma}\label{scen11}
Assume that a smooth solution of \eqref{sqgd1} $\theta(x,t)$ has
modulus of continuity $\omega_B$ at some time $t_0.$ The only way
this modulus of continuity may be violated is if there exists $t_1
\geq t_0$ and $y,z,\,y\not=z$, such that
$\theta(y,t_1)-\theta(z,t_1)=\omega_B(|y-z|),$ while for all $t <
t_1,$ the solution has modulus of continuity $\omega_B.$
\end{lemma}

Next, consider two points $y,z$ and time $t_1$ as in
Lemma~\ref{scen11}. Observe that
\begin{eqnarray}\label{scenest}
\partial_t(\theta(y,t)-\theta(z,t))|_{t=t_1} & =  &u \cdot \nabla
\theta(y,t_1) - u \cdot \nabla \theta(z,t_1)
-(-\Delta)^{1/2}\theta(y,t_1) \\ \nonumber & & +
(-\Delta)^{1/2}\theta(z,t_1)+Au_2(y,t_1)-Au_2(z,t_1).
\end{eqnarray}

Let us denote $|y-z|=\xi.$  We have the following
\begin{lemma}\label{estcru}
For $y,z$ and $t_1$ as in Lemma~\ref{scen11}, we have
\begin{equation}\label{flowest}
|u \cdot \nabla \theta(y,t_1) - u \cdot \nabla \theta(z,t_1)| \leq
\omega_B'(\xi)\Omega_B(\xi)
\end{equation}
and
\begin{equation}\label{dispest}
|Au_2(y,t_1)-Au_2(z,t_1)| \leq A\Omega_B(\xi).
\end{equation}
Moreover, $\delta>\gamma>0$ can be chosen so that
\begin{equation}\label{dissipest}
-(-\Delta)^{1/2}\theta(y,t_1)+ (-\Delta)^{1/2}\theta(z,t_1) \leq
-2\omega_B'(\xi)\Omega_B(\xi).
\end{equation}
\end{lemma}
\begin{proof}
The inequality \eqref{dispest} follows immediately from
Lemma~\ref{somod}. The proof of the inequality \eqref{flowest} is
identical to that provided in \cite{KNV}. The proof of
\eqref{dissipest} is also the same as the treatment of the
dissipative term given in \cite{KNV}. Although the result is not
stated in \cite{KNV} in the same form, it follows immediately from
the arguments provided there. In the estimates above at the point
$x=\delta$ one should use the larger value out of the one-sided
derivatives (which is the left derivative).
\end{proof}

Now we are ready to prove our main technical result, from which
Theorem~\ref{main} follows as explained in the introduction.
\begin{theorem}\label{maint}
Assume that the initial data $\theta_0(x)$ is smooth and periodic.
Then there exists a constant $B(A,\theta_0)$ such that while the
solution of \eqref{sqgd1} $\theta(x,t)$ remains smooth, it
satisfies
\begin{equation}
\|\nabla \theta(\cdot,t)\|_{L^\infty} \leq B.
\end{equation}
\end{theorem}
\begin{proof}
Consider $B_0$ large enough so that $\theta_0(x)$ has $\omega_{B}$
for any $B>B_0.$ Suppose the solution $\theta(x,t)$ loses
$\omega_B,$ then by Lemma~\ref{scen11} we can find $y,z$ and $t_1$
so that $\theta(y,t_1)-\theta(z,t_1)=\omega_B(|y-z|)$ and
$\theta(x,t)$ has $\omega_B$ for all $t \leq t_1.$ By
Lemma~\ref{estcru}, we have
\begin{equation}\label{almfin}
\partial_t(\theta(y,t)-\theta(z,t))|_{t=t_1} \leq
-\omega_B'(\xi)\Omega_B(\xi)+A\Omega_B(\xi).
\end{equation}
Moreover, by Lemma~\ref{Linf}, $\|\theta(\cdot,t)\|_{L^\infty}\leq
D(A,\theta_0)$ and so
\[ \omega_B(\xi) = \theta(y,t)-\theta(z,t) \leq 2D(A,\theta_0). \]
Since $\omega_B(\xi) = \omega(B\xi),$ it follows that $B\xi \leq
\omega^{-1}(2D(A,\theta_0)).$ But then since $\omega'$ is
decreasing,
\[ \omega_B'(\xi) = B \omega'(B\xi) \geq B
\omega'(\omega^{-1}(2D(A,\theta_0)).\] In particular, the right
hand side in \eqref{almfin} is strictly negative if $B \geq
A/\omega'(\omega^{-1}(2D(A,\theta_0)).$ This gives a contradiction
with the definition of $t_1$ since by smoothness the modulus of
continuity should have been violated at an earlier time. Thus
moduli of continuity corresponding to sufficiently large $B$ are
preserved by evolution, as claimed by the Theorem.
\end{proof}

\smallskip

\noindent {\bf Acknowledgement.} \rm Research of AK has been
supported in part by the NSF-DMS grant 0653813. Research of FN has
been partially supported by the NSF-DMS grant 0501067. We thank
Leslie Smith and Jai Sukhatme for suggesting the problem and many
fruitful discussions.

\end{document}